\documentclass[11pt]{amsart}
\usepackage{mathrsfs,url,color,graphicx}
\usepackage{amsthm,eucal,amssymb,amsfonts,amsxtra,amsmath}

\theoremstyle{plain}
    \newtheorem{thm}{Theorem}

    \newtheorem{prop}[thm]{Proposition}

\newcommand{\Fresse}{Fra\"{i}ss\'{e}}

\newtheorem*{claim}{Claim}

\theoremstyle{definition}

\theoremstyle{remark}

\def\ie{i.\,e., }

\newcommand{\ignore}[1]{}

\title{Finite trees are Ramsey \\ under topological embeddings}

\author{Manuel Bodirsky}
   \address{Laboratoire d'Informatique  (LIX), CNRS UMR 7161, \'Ecole Polytechnique}
    \email{bodirsky@lix.polytechnique.fr}
    \urladdr{http://www.lix.polytechnique.fr/~bodirsky}

\author{Diana Piguet}
   \address{Centre for Discrete Mathematics and its Applications (DIMAP)} 
\email{D.Piguet@warwick.ac.uk}
\urladdr{http://www2.warwick.ac.uk/fac/sci/maths/people/staff/Diana\_Piguet}
    \thanks{Diana Piguet was supported by the project 1M0545 by Czech Ministry of Education, and by the FIST (Finite Structures) project, in the framework of the European CommunityÕs \emph{Transfer of Knowledge} programme.}

\begin{document}
\begin{abstract}
We show that the class of finite rooted binary plane trees
is a Ramsey class (with respect to topological embeddings that map leaves to leaves).
That is, for all such trees $P,H$ and every natural number $k$ there exists a tree $T$ such that for every $k$-coloring of the (topological) copies of $P$ in $T$ there exists a (topological) copy $H'$ of $H$ in $T$ such that all copies of $P$ in $H'$ have the same color.
When the trees are represented by the so-called rooted triple relation, the result gives rise to a Ramsey class of relational structures with respect to induced substructures. 
\end{abstract}

   \maketitle

\section{Introduction and Result}
All trees in this paper are finite, rooted, and \emph{binary}, i.e., all vertices have outdegree two or zero. We also assume that the trees are
\emph{plane} in the sense that they are embedded without crossings into the half-plane such that all leaves lie on the boundary of the half-plane (and hence there is a linear order
on the leaves of a tree).
In the following, \emph{tree} always means finite rooted binary plane tree. The set of all vertices of a tree $T$ is denoted by $V(T)$, and the set of all leaves by $L(T)$. 

Two trees $H$ and $T$ are said to be \emph{isomorphic}
if there exists a bijection $f$ from $V(H)$ to $V(T)$ that
preserves the tree structure and that preserves the linear order
on the leaves given by the embedding of the tree.
We say that $H$ is a \emph{(topological) subtree} of $T$ if $L(H) \subseteq L(T)$ and if $T$ can be obtained
from $H$ by adding isolated vertices, adding edges, 
and \emph{subdividing} edges by replacing an edge with a path 
(so that all inner vertices of the path are new vertices).
If $H$ is a subtree of $T$ that is isomorphic to $G$ then
we say that \emph{$H$ is a (topological) copy of $G$ in $T$}.

\begin{figure}
\begin{center}
\includegraphics[scale=0.6]{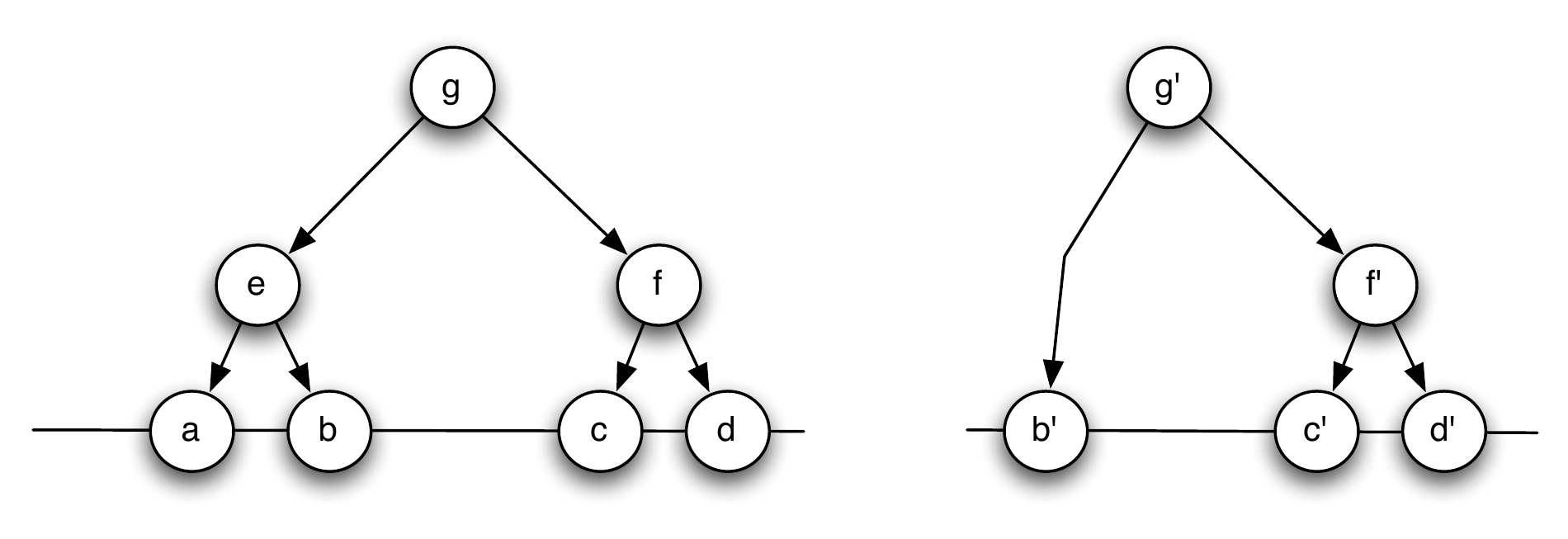} 
\end{center}
\caption{Text}\label{fig:example}
\end{figure}

We illustrate these concepts at the following example, drawn in Figure~\ref{fig:example}. The tree on the left with root $g'$ contains a copy of
the tree on the right: we can subdivide the edge from $g'$ to $b'$
by a new vertex $e'$, add an isolated vertex $a'$, and add an edge from
$e'$ to $a'$. The resulting graph is isomorphic to the graph
on the left. 

A \emph{$k$-coloring} of a set $S$ is a mapping $\chi$ from $S$ into a set of cardinality $k$ (the set of \emph{colors}).
We say that an element of $S$ \emph{has color $c$ (under $\chi$)} if it is mapped to $c$ (by the mapping $\chi$).
In the statement of our result it will be convenient to use the classical
Ramsey theoretic notation $T \rightarrow (H)^P_k$ (for trees $T,P,H$ and an integer $k \geq 2$) for the fact
that for any $k$-coloring of the copies of $P$ in $T$ there is a copy
$H'$ of $H$ in $T$ such that all copies of $P$ in $H'$ have the same color.

We prove the following.

\begin{thm}\label{thm:tramsey}
For all trees $P,H$ and for all $k \geq 2$ there exists a tree $T$
such that $T \rightarrow (H)^P_k$.
\end{thm}

\paragraph{\bf{Note.}}
After this article has been written, we noticed that the main result
is implied by stronger results in~\cite{Milliken} (Theorem 4.3), building
on work in~\cite{DeuberTreeRamsey}. The results in~\cite{Milliken}
also imply results on infinite trees, cover many variations, and therefore
involve more sophisticated terminology than here. Since we
need the result in applications (to study the complexity of 
constraint satisfaction problems), and since we found both the statement
and its proof easier accessible in the present paper, 
we decided to still post this article as a technical report,
without claiming originality.

\subsection{Ramsey classes.}
The result is of interest in structural Ramsey theory since it gives rise to a so-called \emph{Ramsey class} of relational structures with respect to embeddings. Ramsey classes of relational structures
are one of the central topics in Ramsey theory, and we give a brief introduction.

A \emph{relational signature} $\tau$ is a set of relation symbols $R$, each associated with an arity $\textit{ar}(R) \geq 1$. 
A \emph{relational structure} $\Gamma$ with signature $\tau$ (short, a $\tau$-structure) consists
of a \emph{domain} $D_\Gamma$ and a relation $R^\Gamma \subseteq (D_\Gamma)^k$ for each relation symbol $R \in \tau$ of arity $k$.

Let $\Gamma$ and $\Delta$ be two relational structures over the same
signature $\tau$. Then an \emph{embedding} of $\Gamma$ into $\Delta$ is an injective mapping $f$ from $D_\Gamma$ to $D_\Delta$
that satisfies that  $(t_1,\dots,t_{{\textit{ar}(R)}}) \in R^\Gamma$ if and only if
$(f(t_1),\dots,f(t_{{\textit{ar}(R)}})) \in R^\Delta$, for each relation symbol $R \in \tau$. A bijective embedding is called an \emph{isomorphism}.
All classes of structures in this paper will be closed under taking
isomorphisms. If $f$ is the identity mapping, then $\Gamma$ is called
an \emph{(induced) substructure} of $\Delta$. 

A class of relational $\tau$-structures $\mathcal C$ is a \emph{Ramsey class} (with respect to embeddings) if for all $P,H \in {\mathcal C}$ and every natural number $k$ there exists a $T \in \mathcal C$ such that
for every $k$-coloring of the substructures of $T$ that are isomorphic to $P$ there exists a substructure $H'$ of $T$ that is isomorphic to $H$
such that all substructures of $H'$ that are isomorphic to $P$ have the same color; again, this is denoted by  $T \rightarrow (H)^P_k$.
Examples of Ramsey classes are 
\begin{itemize}
\item the class of all finite structures over the empty signature (this is just the classical theorem of Ramsey);
\item the class of all finite linear orders over the signature $\{<\}$ with a single binary relation symbol $<$ that defines the ordering (yet another form of the classical theorem of Ramsey);
\item for any relational signature $\tau$, the class of all ordered structures over $\tau$ (i.e., all finite structures over the signature $\tau \cup \{<\}$, where $<$ is a new binary relation symbol that defines the ordering)~\cite{NesetrilRoedlPartite,AbramsonHarrington,NesetrilRoedlOrderedStructures};
\item the class of all finite posets that are equipped with a linear extension~\cite{PosetRamsey};
\item the class of all finite metric spaces~\cite{MetricRamsey}.
\item Canonically ordered Boolean Algebras~\cite{Topo-Dynamics}.
\end{itemize}
For more examples of Ramsey classes, 
see e.g.~\cite{RamseyClasses,NesetrilSurvey,Topo-Dynamics}. It is known that all Ramsey classes $\mathcal C$ of finite structures that are closed
under taking induced substructures (those classes are sometimes also called \emph{hereditary}) are \emph{amalgamation classes}~\cite{RamseyClasses}, and hence there
exists an (up to isomorphism unique) \emph{homogeneous}Ê
countably infinite structure $\Gamma$ such that $\mathcal C$ is exactly the class of finite structures that
embed into $\Gamma$ (this follows from \Fresse's theorem, see~\cite{Hodges}).
Homogeneity is a very strong model-theoretic property, 
and it is therefore possible to use model-theoretic techniques and results to approach a classification of all Ramsey classes 
that are closed under taking induced substructures.
This is the program that has basically been launched in~\cite{RamseyClasses}.

The result presented here gives rise to a new Ramsey class, and hence contributes to the classification program. We have to describe how to represent rooted binary plane trees
as relational structures. Our relational structures will be ordered,
i.e., the signature contains a binary relation symbol $<$ that is interpreted by a linear order. The tree structure is represented
by a single ternary relation symbol as follows.
For leaves $a,b,c$ of a tree, we write $ab|c$ if the least common ancestor of $a$ and $b$ is below the least common ancestor of $a$ and $c$ in the tree; the relation $|$ is also called the \emph{rooted triple relation}, following terminology in phylogenetic reconstruction~\cite{Steel,NgWormald,HenzingerKingWarnow}. It is known that a rooted binary tree is described up to isomorphism by the rooted triple relation (see e.g.~\cite{Steel}).

We now associate to a rooted binary plane tree $T$ the
relational structure $\Gamma = (L(T); |,<)$ where a triple $(a,b,c)$ 
of elements of $\Gamma$ is
in the rooted triple relation $|^\Gamma$ if $T$ satisfies $ab|c$,
and where a pair of elements $(a,b)$ is in the relation $<^\Gamma$
if the leaf $a$ lies on the left of the leaf $b$ with respect to the embedding of $T$ into the half-plane.
It is now clear that Theorem~\ref{thm:tramsey} implies that the class of relational structures $\mathcal T$ over the signature $\{|,<\}$ obtained as described above from rooted binary plane trees is a Ramsey class.

\subsection{Infinite Permutation Groups}
As we have mentioned in the previous subsection, any Ramsey class that is closed under taking substructures (and our Ramsey class $\mathcal T$ is obviously closed under taking substructures) is an amalgamation class, and therefore
there exists a unique countable homogeneous structure $\Lambda = (\Omega; |, <)$ such that $\mathcal T$ is exactly the class of all finite structures that embed into $\Lambda$. The structure $(\Omega; |)$ 
(i.e., the reduct of $\Lambda$ that only contains the rooted triple relation without an ordering on the domain) is well-known to model-theorists and in the theory of infinite permutation groups, and also has many explicit constructions; see e.g.~\cite{Oligo,AdelekeNeumann}.
Its automorphism group is oligomorphic, 2-transitive, and 3-set-transitive, but not 3-transitive. The rooted triple relation $|$ is 
a $C$-relation in the terminology of~\cite{AdelekeNeumann}.

\section{Proof of the main result}
We start with the easy
special case where we color only the leaves of a tree (and thus the copies of trees of order $1$); this will serve us as an induction basis in the proof of the main result.
We denote by $T(c)$ the rooted binary tree with $2^c$ leaves of height $c$, \ie any leaf is at distance $c$ to the root. 
In particular, $T(0)$ denotes the one-vertex tree.

\begin{prop}\label{prop:ramsey-singelton}
For all trees $H$ and all $k \geq 2$ there exists a tree $T$
such that $T \rightarrow (H)^{T(0)}_k$.
\end{prop}

\begin{proof}
Note that the $k$-colorings of the copies of $T(0)$ in a tree $T$
are just colorings of the leaves of $T$; we therefore just speak
of $k$-colorings of $T$.

We apply the following operation to construct trees.
Let $G$ and $H$ be trees. Then $G[H]$ denotes the tree
obtained from $G$ by replacing each leaf of $G$ by a tree isomorphic to  $H$ (for each leaf $v$ of $G$ the children of the root of the copy of $H$
in $G[H]$ are the children of $v$). It is clear that the resulting
tree has an embedding into the half-plane so that $G[H]$ is again a rooted binary plane tree.
We can iterate the construction: let $H^{(1)}$ be $H$,
and define $H^{(i+1)}$ for $i \geq 1$ to be $H[H^{(i)}]$.

Clearly, $H[H] \rightarrow (P)^{T(0)}_2$, because for all $2$-colorings of $H[H]$ either
one of the `lower' copies of $H$ in $H[H]$ (i.e., one of the copies 
of $H$ in $H[H]$ that replaced a leaf in $H$) is monochromatic, 
and we are done, or otherwise all these copies of $H$ contain both colors.
Let $a_1,\dots,a_n$ be the leaves
of $H$, and let $c_i$ be a leaf in the copy
of $H$ that replaced $a_i$ in $H[H]$ and that has color $0$. 
Then the subtree of $H[H]$ with leaf set $\{c_1, \dots, c_n\}$ is a 
$0$-chromatic copy of $H$.

Therefore
\begin{align*}
H[H^{(2)}] & \rightarrow (H)^{T(0)}_3
\end{align*}
because either one of the lower copies of $H^{(2)}$ in $H[H^{(2)}]$
is 2-chromatic, in which case we have already shown that
this copy contains a monochromatic copy of $H$, or all copies contain
all three colors. But then, by an analogous argument as above, we have a monochromatic copy of $H$ in $H^{(3)}$. By iterating this argument it follows
that $H^{(k)} \rightarrow (H)^{T(0)}_k$.
\end{proof}

In the proof of Theorem~\ref{thm:tramsey}, the following notation will be convenient.
The set of all copies of $P$ in $H$ is denoted by $H \choose P$.
If all copies of $P$ in $H$ have the same color, we say that $H$ is \emph{$\chi$-monochromatic} (or simply monochromatic if the coloring is clear from the context). If the color is $k$, we also say that $H$ is \emph{$k$-chromatic}.

If $T$ is a tree with more than one vertex, then the root of $T$ has
exactly two children; we denote the subtree $T$ rooted at the left child
by $T_{\swarrow}$, and the subtree of $T$ rooted at the right child
by $T_{\searrow}$ (and we speak of the \emph{left subtree of $T$} and the \emph{right subtree of $T$}, respectively).
Finally, suppose that $H_1$ and $H_2$ are disjoint subtrees of $T$. Then
$\left<H_1,H_2 \right>$ denotes the (uniquely defined) subtree of $T$ with leaves
$L(H_1) \cup L(H_2)$.

\begin{proof}[Proof of Theorem~\ref{thm:tramsey}]
We prove the theorem by induction on the size of $P$. For $P=T(0)$ the statement holds by Proposition~\ref{prop:ramsey-singelton}. We now assume that the statement holds for all proper subtrees of $P$; we want to prove it for $P$. 

We start with the case $k=2$, and proceed by induction on the size of $H$. Observe that trivially $P\rightarrow (P)^P_2$. 
So assume that the theorem holds for proper subtrees $H'$ of $H$, \ie
we assume that there exists a tree $T_{H'}$ such that
 $T_{H'} \rightarrow (H')^P_2$.
In particular, we assume that the theorem holds for $H_\swarrow$ and $H_\searrow$, the left and right subtree of $H$.
In the inductive step, we first prove the following claim.

\begin{claim}[Asymmetric step]\label{claim:asym}
There exists a tree $F$ such that for any $2$-coloring $\chi: {F \choose P} \rightarrow \{0,1\}$
of the copies of $P$ in $F$
\begin{itemize} 
\item there is a $0$-chromatic copy of $H_\searrow$ in $F_\searrow$ or of $H_\swarrow$ in $F_\swarrow$, or
\item there exists a 1-chromatic copy of $H$ in $F$.
\end{itemize}
\end{claim}

\begin{proof}[Proof of the asymmetric step]
For a sufficiently large $n$ (whose choice will be discussed at the end of the proof), let $F$ be isomorphic to $T(n)$. 
Suppose that there is no $0$-chromatic copy of $H_\swarrow$ and no $0$-chromatic copy of $H_\searrow$ in $F$ under the coloring $\chi$. We show that there exists a $1$-chromatic copy of $H$ in $F$. 
Let $\psi:{F_\swarrow\choose P_\swarrow}\rightarrow 2^{F_\searrow \choose P_\searrow}$ be the mapping 
 that assigns to each copy of $P_\swarrow$ in $F_\swarrow$ the coloring of ${F_\searrow \choose P_\searrow}$ induced in the following way:
if $P_1$ is a copy of $P_\swarrow$ in $F_\swarrow$, then
$\psi(P_1)$ is the mapping that maps a copy $P_2$ of $P_\searrow$ in $F_\searrow$ to $\chi(\left< P_1,P_2 \right>)$.
Hence, we color the set 
${F_\swarrow\choose P_\swarrow}$ where the colors are themselfes 2-colorings of ${F_\searrow \choose P_\searrow}$. 
 
By inductive hypothesis, and because $F_\swarrow$ is large enough, 
there exists a $\psi$-monochromatic copy $F_1$ of $T(m)$ in $F_\swarrow$, where $m$ is sufficiently large. 
Let $\phi$ be the color of the copies of $P_\swarrow$ in $F_1$; recall
that $\phi$ is a $2$-coloring of ${F_\searrow \choose P_\searrow}$. Since 
$F_\searrow$ is large enough, there is a copy $F_2$ of $T(m)$ that is $\phi$-monochromatic, 
say all copies of $P_\searrow$ in $F_2$ are \textit{blue}. 
Note that if $P_1$ is a copy of $P_\swarrow$ in $F_1$,
and $P_2$ is a copy of $P_\searrow$ in $F_2$, then
$\left < P_1,P_2 \right>$ is a copy of $P$ in $F$ that is colored \textit{blue}.

Assume first that $\textit{blue} = 1$. By inductive assumption, $F_1 \rightarrow (H_\swarrow)^P_2$ and 
$F_2 \rightarrow (H_\searrow)^{P}_2$. The color of the copies of $P$ in the monochromatic copy $H_1$ of $H_\swarrow$ in $F_1$ and in the monochromatic copy $H_2$ of $H_\searrow$ in $F_2$ must be $1$, or otherwise the first disjunct of the conclusion of the statement is fulfilled. Because $\left< H_1,H_2 \right>$ is a copy of $H$ in $F$, all copies of $P$ in $\left< H_1,H_2 \right>$ 
are colored by $1$, and we are done in this case.
So we can assume that $\textit{blue} = 0$.

We now iterate this argument $h$-times as follows, where
$h$ is the height of $H$ (the maximal distance from the root of $H$ to one of its leaves).
In the $i$-th step, we define disjoint subtrees $F_1^i, \dots, F^i_{2^i}$ of $F$. Initially, in the first step, let $F^1_1 := F_1$ and $F^1_2 := F_2$.
In the $i+1$-st step, for $i \geq 1$ and $j \leq 2^i$, let $\psi^i_j$ 
be the mapping that assigns to each copy of $P_\swarrow$ in $(F^i_j)_\swarrow$ the coloring of ${(F^i_j)_\searrow \choose P_\searrow}$ induced by $\chi$ as before.
  
By inductive hypothesis, and because $(F^i_j)_\swarrow$ is large enough, 
there exists a $\psi^i_j$-monochromatic copy $F^{i+1}_j$ of $T(n_i)$ in $(F^{i}_j)_\swarrow$, where $n_i$ is sufficiently large. 
Let $\phi^i_j$ be the color of the copies of $P_\swarrow$ in $F^{i+1}_j$. Since $(F^{i}_j)_\searrow$ is large enough, there is a copy $F^{i+1}_{j+2^i}$ of $T(n_i)$ that is $\phi^i_j$-monochromatic.
We can argue as before to conclude that 
if $P_1$ is a copy of $P_\swarrow$ in $F^{i+1}_j$
and $P_2$ is a copy of $P_\searrow$ in $F^{i+1}_{j+2^i}$, then
$\left < P_1,P_2 \right>$ is a copy of $P$ in $F$ that has color $0$.
Finally, we select one leaf in each of the trees $F^h_1, \dots,F^h_{2^h}$.  
These vertices show that there is a copy of $T(h)$ and hence also a copy of $H$ in $F$ where every copy of $P$
is colored by $0$.

We can certainly find appropriate (large) values for $n,m$, and $n_i$, for $i \geq 1$, since we can choose $n_h = 1$, and for $i < h$ we can choose $n_i$ large enough depending
on the size of $n_{i+1}$, so that we can finally also 
choose an appropriate value for $m$ and for $n$.
\end{proof}

To conclude the inductive proof of Theorem~\ref{thm:tramsey}, let $T$ be a copy of $T(d)$
 where $d$ is large enough (again we discuss the choice of $d$ at the end of the proof).
We will show that for any $\chi: {T \choose P}\rightarrow \{0,1\}$ there exists a monochromatic copy of $H$. Let $\psi:{T_\swarrow\choose P_\swarrow}\rightarrow 2^{T_\searrow \choose P_\searrow}$ be the function that assigns 
to a copy $P_1$ of $P_\swarrow$ in $T$ the function that maps a copy $P_2$ of $P_\searrow$ in $T_\searrow$ to $\chi(\left< P_1,P_2 \right>)$.
By our inductive hypothesis, and since $T_\swarrow$ is large enough, we find a $\psi$-monochromatic copy $T_1$ of the tree $F$ given by the assymetric step.
This gives us a 2-coloring $\phi$ of ${T_\searrow \choose P_\searrow}$. Since $T_\searrow$ is large enough, we find a $\phi$-monochromatic 
copy $T_2$ of the tree $F$ from the asymmetric step;
let us assume that all copies of $P$ in $T_2$ are colored by $1$.
Note that if $P_1$ is a copy of $P_\swarrow$ in $T_1$ and
$P_2$ is a copy of $P_\searrow$ in $T_2$, then $\left< P_1, P_2\right>$ is a subtree of $T$ that is colored by $1$ under $\chi$.

We apply the asymmetric step 
to $T_1$ and $T_2$ and the restriction of $\chi$ to $T_1$ and $T_2$, respectively. If there is a $1$-chromatic copy of $H$ in $T_1$ or in $T_2$, we are done. So we may assume that the colorings of ${T_1\choose P}$ and ${T_2 \choose P}$ are such that there is a $0$-chromatic copy $H_1$ of $H_\swarrow$ in $T_1$ and a $0$-chromatic copy $H_2$ of $H_\searrow$ in $T_2$. Then $\left< H_1, H_2 \right>$ is a copy of $H$ in $T$ and $1$-chromatic with respect to $\chi$.

We have proved that there exists a tree $T$ such that 
$T \rightarrow (H)^P_2$, and now prove the theorem for any finite number of colors $k$.
Let $l$ be $\lceil \log_2 k\rceil$.
Define $T^0$ to be $H$, and let $T^i$ for $1 \leq i \leq l$
be such that $T^i \rightarrow (T^{i-1})_2^P$; we already know that such a tree $T^i$ exists. We claim that $T^l \rightarrow (H)_k^P$.
Let $\chi: {T^l \choose P} \rightarrow \{0,\dots,k-1\}$, and consider for $1 \leq i \leq l$ the colorings $\psi_i: {T^i \choose P} \rightarrow \{0,1\}$
where $\psi_i$ colors a copy $P'$ of $P$ by $1$ if the $i$-th bit in the binary representation of $\chi(P')$ is 1, and it colors $P'$ by $0$ otherwise.
For $1 \leq i \leq \lceil \log_2 k\rceil$, let $S_{i-1}$ be the $\psi_i$-monochromatic copy of $T^{i-1}$ in $S_i$, and let $b_i$ be the
color of the copies of $P$ in $S_{i-1}$. Note that $S_0$ is isomorphic to $H$. All copies of $P$ in $S_0$ have color $b_i$ with respect to $\psi_i$, for all $0 \leq i \leq \lceil \log_2 k\rceil$, and
hence they have the same color with respect to $\chi$.
\end{proof}



\paragraph{\bf Remark.} With minor modifications of the proof
a similar result can be shown for the class of trees with respect
to embeddings where leaves are not necessarily mapped to leaves.

\bibliographystyle{alpha}
\bibliography{../local}

We appologize for the improper spelling of the name of Jaroslav Ne\v{s}et\v{r}il in three of the references; the latex system does not allow for the symbol \v{s} in the automatically generated short-cuts for the citations.

\end{document}